\documentclass[12pt,draftcls,onecolumn]{IEEEtran}
\IEEEoverridecommandlockouts                              
\usepackage{color}
\usepackage{graphicx}
\usepackage{subfigure}
\usepackage{amsfonts}
\usepackage{amsmath}
\usepackage{amssymb}
\usepackage{psfig}

\newcommand \nn \nonumber
\newcommand \nnd \noindent

\newtheorem{thm}{Theorem}

\title{Coverage Optimization using Generalized Voronoi Partition}

\author{K.R. Guruprasad
\thanks{
K.R. Guruprasad is an Assistant Professor at the Department of
Mechanical Engineering, National Institute of Technology Karnataka,
Surathkal, 575025, India, and is currently with the Department of
Computer Science, University of Nebraska, Omaha, Nebraska, USA, as a
visiting faculty. {\tt\small krgprao@gmail.com} ({\bf Corresponding
author})} and Debasish Ghose
\thanks{Debasish Ghose
is a Professor at the Guidance, Control, and Decision Systems
Laboratory, Department of Aerospace Engineering, Indian Institute of
Science, Bangalore, 560012, India. {\tt\small
dghose@aero.iisc.ernet.in}}}

\begin{document}

\maketitle


\begin{abstract}
In this paper a  generalization of the Voronoi partition is used for
optimal deployment of autonomous agents carrying sensors with
heterogeneous capabilities, to maximize the sensor coverage. The
generalized centroidal Voronoi configuration, in which the agents
are located at the centroids of the corresponding generalized
Voronoi cells, is shown to be a local optimal configuration.
Simulation results are presented to illustrate the presented
deployment strategy.
\end{abstract}


\section{Introduction}
Technological advances in areas such as wireless communication,
autonomous vehicular technology, computation, and sensors,
facilitate the use of a large number of agents (UAVs, mobile robots,
autonomous vehicles etc.) to cooperatively achieve various tasks in
a distributed manner. One of the very useful applications of the
multi-agent systems is in sensor networks, where a group of
autonomous agents perform cooperative sensing of a large
geographical area. Voronoi partition \cite{okabe} has been used for
optimal deployment of mobile agents carrying sensors
\cite{bullo1,tase,ijss1}. This paper discusses a generalization of
Voronoi partition to addresses heterogeneity of sensors in their
capabilities.

The problem of optimal deployment of sensors has attracted many
researchers \cite{li2}. Zou and Chakrabarty \cite{zou} use the
concept of virtual force to solve this problem.  Akbarzadeh et al.
\cite{akbarzadeh} use an evolutionary algorithm approach. Yao et al.
\cite{yao} address a problem of camera placement for persistent
surveillance. Murray et al. \cite{atmurray} address a problem of
placing sensors such as video cameras for security monitoring. Kale
and Salapaka \cite{Kale11} use heuristically designed algorithms
based on maximum entropy to seek global minimum for a simultaneous
resource location and multi-hop routing. The problem of optimal
deployment of sensors belongs to a class of problems known as
locational optimization or facility location \cite{okabe,drezner} in
the literature, and when homogeneous sensors are used, a centroidal
Voronoi configuration \cite{du} is a standard solution for this
class of problems.  Cortes et al. \cite{bullo1,bullo2} use these
concepts to formulate and solve the problem of optimal deployment of
sensors. The authors provide rigorous mathematical results on
spatial distribution, convergence, and other useful properties.
Pimenta et al. \cite{vkumar} follow a similar approach to address a
problem with heterogeneous robots, using power diagram (or Voronoi
diagram in Laguerre geometry) to account for different footprints of
the sensors. Kwok and Mart\'{i}nez \cite{kwok} use power weighted
and multiplicatively weighted Voronoi partitions to solve an energy
aware limited range coverage problem. Pavone et al. \cite{pavone}
use power diagrams for equitable partitioning for robotic networks.



In the literature, sensors in a network are usually considered to be
homogeneous in their capability. Whereas, in practical problems, the
sensors may have different capabilities even though they are similar
in their functionality. The heterogeneity in capabilities could be
due to various reasons, the chief being the difference in specified
performance. Though some researchers address optimal deployment of
heterogeneous sensors (such as \cite{vkumar}), one of the well known
generalization of the standard Voronoi partition (weighted Voronoi
partition, Power diagram, etc.) is used to address heterogeneity. In
this paper, we generalize the standard Voronoi partition to address
a class of heterogeneous locational optimization problems, which
include some of the problems addressed in the literature. We use
node functions in place of the usual distance measure used in
standard Voronoi partition and its variants. The mobile sensors are
assumed to have heterogeneous capabilities in terms of their
effectiveness. This paper presents a generalization of the optimal
deployment concepts presented in \cite{bullo1,bullo2}, using
generalized Voronoi partition in place of the standard Voronoi
partition. Some preliminary results have been reported in
\cite{isvd08} earlier. This paper gives a more elaborate treatment.

%

\section{Generalization of the Voronoi partition}
A generalization of the Voronoi partition, considering heterogeneity
in the sensors' capabilities, is presented here. Several extensions
or generalizations of Voronoi partition to suit specific
applications have been reported in the literature \cite{okabe}.
Herbert and Seidel \cite{herbert}  have introduced an approach in
which, instead of the site set, a finite set of real-valued
functions $f_i: D \mapsto \mathbb{R}$ were used to partition the
domain $D$.

Let $I_N = \{1,2,\ldots,N\}$. Consider a space $Q \subset
\mathbb{R}^d$ and a set of points called {\em nodes} or {\em
generators} $\mathcal{P} = \{p_1,p_2, \ldots, p_N \}$, $p_i \in Q$,
with $p_i \neq p_j$, whenever $i \neq j$. The {\em Voronoi
partition} generated by $\mathcal{P}$ is the collection
$\mathcal{V}(\mathcal{P}) = \{V_i(\mathcal{P})\}_{i \in I_N}$, and
is defined as,
\begin{equation}
\label{vor_def} V_i(\mathcal{P}) = \left \{ q \in Q | \parallel
q-p_i
\parallel \leq
\parallel q - p_j \parallel, \forall p_j \in \mathcal{P}, j \in I_N \right\}
\end{equation}
where, $\parallel . \parallel$ denotes the Euclidean norm.
%
%
The basic components of the Voronoi partition are: i) A space to be partitioned; ii) A set of sites, or nodes, or generators;
and iii) A distance measure such as the Euclidean distance.


In this paper a generalization of the Voronoi partition is defined
to suit the application, namely the heterogeneous locational
optimization of sensors. This generalization has the following
components: i) The domain of interest as the space to be
partitioned; ii) The configuration of multi-agent system
$\mathcal{P}$ as the site (or node or generator) set; and iii) A set
of node functions in place of a distance measure. Consider strictly
decreasing analytic functions $f_i : \mathbb{R}^+ \mapsto
\mathbb{R}$, where, $f_i$ is called a {\em node function} for the
$i$-th node. Define {\em generalized Voronoi partition} of $Q$ with
node configuration $\mathcal{P}$ and node functions $f_i$ as a
collection $\mathcal{V}^g = \{V^g_i\}$, $i\in I_N$, with mutually
disjoint interiors, such that $Q=\cup_{i\in I_N} V^g_i$, where
$V^g_i$ is defined as
\begin{equation}
\label{vor_fun} V^g_i = \{ q\in Q  |  f_i(\| p_i - q \|) \geq f_j(\|
p_j - q \|), \quad \forall j \neq i\text{,} j \in I_N \}
\end{equation}
A generalized Voronoi cell $V^g_i$ can be topologically
non-connected, null, and may contain other Voronoi cells. In the
context of the problem discussed in this paper, $q \in V^g_i$ means
that the $i$-th agent/sensor is the most effective in sensing at
point $q$. This is reflected in the $\geq$ sign in the definition.
In standard Voronoi partition used for the homogeneous case, $\leq$
sign for distances ensured that the $i$-th sensor is the most
effective in $V_i$. (Note that $f_i$ are strictly decreasing.) The
condition that $f_i$ are analytic implies that for every $i,j \in
I_N$, $f_i - f_j$ is analytic. By the property of real analytic
functions, the set of intersection points between any two node
functions is a set of measure zero. This ensures that the
intersection of any two cells is a set of measure zero, that is, the
boundary of a cell is made up of the union of at most $d-1$
dimensional subsets of $\mathbb{R}^d$. Otherwise the requirement
that the cells should have mutually disjoint interiors may be
violated. Analyticity of the node functions $f_i$ is a sufficient
condition to discount this possibility. Further, note that for some
$i \neq j$, if $f_i \neq f_j$, then it is acceptable to have $p_i =
p_j$. Only when $f_i = f_j$ and $i \neq j$, the condition $p_i \neq
p_j$ is imposed in case of the generalized Voronoi partition (see
Appendix for details). It can be shown that the centroid of a
generalized Voronoi cell lies within its convex hull \cite{Bullo09}.

\subsection{Special cases~}
A few interesting special cases of the generalized Voronoi partition
are  discussed below.

\subsubsection*{Case 1: Weighted Voronoi partition} Consider
multiplicatively and additively weighted Voronoi partitions as
special cases. Let $f_i(r_i) = -\alpha_i r_i - d_i$
where, $r_i = \parallel p_i-q \parallel$ and,  $\alpha_i$ and
$d_i$ take finite positive real values for $i \in I_N$.
Thus,
\begin{equation}
\label{gen_mult_vor} V^w_i = \{ q \in Q  |  \alpha_i r_i + d_i \leq \alpha_j r_j + d_j \text{,}\quad \forall j \neq i \text{,} \quad j
\in I_N \}
\end{equation}
The partition $\mathcal{V}^w = \{V^w_i\}$ is called a
multiplicatively and additively weighted Voronoi partition.
$\alpha_i$ are called multiplicative weights and $d_i$ are called
additive weights.

\subsubsection*{Case 2: Standard Voronoi partition} The standard
Voronoi partition can be obtained as a special case of
(\ref{vor_fun}) when the node functions are  $f_i(r_i)= -r_i$. It can be shown that if the node functions are homogeneous
($f_i(\cdot) = f(\cdot)$ for each $i \in I_N$), then the
generalized Voronoi partition gives the standard Voronoi
partition.

\subsubsection*{Case 3: Power diagram}Power diagram or Voronoi diagram in Laguerre geometry
$\mathcal{V}^{LG} = \{V^{LG}_i\}$ is defined as $V^{LG}_i = \{q \in
Q | d_p(q,p_i) \leq d_p(q,p_j)\text{,} i\neq j \}, i,j \in I_N$
where, $d_p(q,p_i) = \|p_i-q\|^2 - R^2_{p_i}$, the power distance
between $q$ and $p_i$, with $R_{p_i}>0$ being a parameter fixed for
a given node $p_i$. In the context of robot coverage problem
addressed in \cite{vkumar}, $R_{P_i}$ represents the radius of the
footprint of the $i$-th robot. It is easy to see that the power
diagram can be obtained from the generalized Voronoi partition
(\ref{vor_fun}) by setting $f_i(\|q-p_i\|) = -(\|p_i-q\|^2 -
R^2_{p_i})$ with $R_{P_i}$ as a parameter specific to each node
function.

%

\subsubsection*{Case 5: Other possible variations}
Other possible variations of the Voronoi partition are using objects
such as lines, curves, discs, polytopes, etc. other than points as
sites/nodes, generalization of the space to be partitioned, and use
of non-Euclidean metrics or pseudo-metrics. It is easy to verify
that these generalizations can also be obtained by suitable
selection of site sets, spaces, and node functions.


\section{Heterogeneous locational optimization problem}
The heterogeneous locational optimization problem (HLOP) for a
mobile sensor network is formulated and solved here. Let $Q \subset
\mathbb{R}^d$ be a convex polytope, the space in which the sensors
have to be deployed; $\phi :Q \mapsto [0,1]$, be a density
distribution function, with $\phi(q)$ indicating the probability of
an event of interest occurring at $q\in Q$, indicating the
importance of measurement at $q$; $\mathcal{P}(t) = \{p_1(t),p_2(t),
\ldots, p_N(t)\}, i \in I_N$, $p_i(t) \in Q$ be the configuration of
$N$ sensors at time $t$, with $p_i(0) \neq p_j(0), \forall i\neq j$;
$f_i:\mathbb{R}^+ \mapsto \mathbb{R}$, $i \in I_N$, be analytic,
strictly decreasing function corresponding to the $i$-th node, the
sensor effectiveness function of $i$-th agent, with $f_i(\|p_i-q\|)$
indicating the effectiveness of the $i$-th sensor located at $p_i
\in Q$, in sensing at a point $q \in Q$. It is natural to assume
$f_i$ to be strictly decreasing. The objective of the problem is to
deploy the mobile sensors in $Q$ so as to maximize the probability
of detection of an event of interest.

In case of homogeneous sensors \cite{bullo1}, the sensor located in
Voronoi cell $V_i$ is closest to all the points $q\in V_i$ and
hence, by the strictly decreasing variation of sensor's
effectiveness with distance, most effective within $V_i$. Thus, the
Voronoi decomposition leads to optimal partitioning of the space in
the sense that each sensor is most effective within its Voronoi
cell. In the heterogeneous case too, it is easy to see that each
sensor is most effective within its generalized Voronoi cell. Now,
as the partitioning is optimal, the problem is to find the location
of each sensor within its generalized Voronoi cell.

\subsection*{The objective function}
Consider the following objective function to be maximized,
\begin{equation}
\label{HLOP} \mathcal{H}(\mathcal{P}) = \int_Q
\max_{i \in I_N}\{f_i(\|q-p_i\|)\}\phi(q)dQ = \sum_{i \in I_N}\int_{V^g_i}f_i(\|q-p_i\|)\phi(q)dQ
\end{equation}
Note that the generalized Voronoi decomposition splits the objective
function into a sum of contributions from each generalized Voronoi
cell. Hence the optimization problem can be solved in a spatially
distributed manner, that is, the optimal configuration can be
achieved by each sensor solving only that part of the objective
function which corresponds to its own cell, using only local
information. Note also that the objective function has the same form
as in \cite{bullo1} where a similar problem is addressed for
homogeneous sensors.

By generalizing results in \cite{Bullo09} (Proposition 2.23 in
\cite{Bullo09}), it can be shown that the gradient of the objective
function (\ref{HLOP}) with respect to $p_i$ is
\begin{equation}
\label{grad} \frac{\partial \mathcal{H}}{\partial p_i} =
\int_{V^g_i}\phi(q)\frac{\partial f_i(r_i)}{\partial p_i}dQ
\end{equation}

\subsection*{The critical points}
By suitably generalizing the concepts in \cite{Bullo09}, it can be
shown that the gradient of the objective function (\ref{HLOP}) with
respect to $p_i$  is
\begin{equation}
\label{grad_HLOP} \frac{\partial \mathcal{H}}{\partial p_i} =
\int_{V^g_i}\phi(q)\frac{\partial f_i(r_i)}{\partial p_i}dQ
= -\int_{V^g_i}\tilde{\phi^i}(q,p_i)(p_i-q)dQ = -\tilde{M}_{V^g_i}(p_i -
\tilde{C}_{V^g_i})
\end{equation}
where, $\tilde{\phi^i}(q,p_i) = -2\phi(q)\frac{\partial
f_i(r_i)}{\partial {(r_i}^2)}$. As $f_i$s are strictly decreasing,
$\tilde{\phi}(q)$ is always non-negative. Hence, $\tilde{\phi}$ can
be interpreted as density modified or perceived by the sensors,
$\tilde{M}_{V^g_i}$ as mass, and $\tilde{C}_{V^g_i}$ as centroid of
the cell $V^g_i$. Thus, the critical points/configurations are
$\mathcal{P}^* = \{\mathcal{P} | p_i =
\tilde{C}_{V^g_i}(\mathcal{P}), i \in I_N\}$, and such a
configuration $\mathcal{P}^*$, of agents is called a {\em
generalized centroidal Voronoi configuration}.



The critical points are not unique, as with the standard Voronoi
partition. But in the case of the generalized Voronoi partition,
some of the cells could become null. Further, as in the case of
homogenous sensors (that is, using standard Voronoi partition), the
critical points/configuration correspond to local minima of the
objective function (\ref{HLOP}).

\subsection{The control law}
Consider the agent dynamics and control law as
\begin{eqnarray}
\label{dyn1_HLOP} \dot p_{i}(t) &=& u_i(t)\\
\label{ctrl1_HLOP} u_i(t) &=& -k_{prop}(p_{i}(t) - \tilde{C}_{V_{i}}(t))
\end{eqnarray}Control law (\ref{ctrl1_HLOP}) makes the mobile sensors move
toward $\tilde{C}_{V^g_{i}}$ for positive $k_{prop}$. If, for some
$i \in I_N$, $V^g_i = \emptyset$, then set $\tilde{C}_{V^g_i} =
p_i$. We restrict the analysis to the simple first order dynamics
under the assumption, that there is a low-level controller which can
cancel the actual dynamics of the agents and enforce
(\ref{dyn1_HLOP}).
%
%

%

\begin{thm}
\label{LaSalle_HLOP} The trajectories of the sensors governed by
the control law (\ref{ctrl1_HLOP}), starting from any initial
condition $\mathcal{P}(0) $, will asymptotically converge to the
critical points of $\mathcal{H}$.
\end{thm}

\noindent {\it Proof.~} Consider $V(\mathcal{P}) = -\mathcal{H}$.
\begin{equation}
\label{vdot} {\dot V}(\mathcal{P}) =
-\frac{d\mathcal{H}}{dt}
=  -\sum_{i \in I_N} \frac{\partial \mathcal{H}}{\partial
p_{i}}\dot{p}_{i}
= -k_{prop}\sum_{i \in I_N}  \tilde{M}_{V^g_{i}}(p_{i} -
\tilde{C}_{V^g_{i}})^2
\end{equation}
Let $p_i \in \partial Q$, the boundary of $Q$ for some $i$ at some
time $t$. The vector $\dot p_i = -k_{prop}(p_i -
\tilde{C}_{V^{g}_i})$ always points inward to $Q$ or is tangential
to $\partial Q$ at $p_i$ as $\tilde{C}_{V^g_i} \in Q$ (Note that
$\text{co}(V^g_i) \subset Q$). Thus, by Theorem 3.1 in
\cite{blanchini}, the set $Q$ is invariant under the closed-loop
dynamics given by (\ref{dyn1_HLOP}) and (\ref{ctrl1_HLOP}).

Further, observe that,  $V: Q^N\mapsto \mathbb{R}$ is continuously
differentiable in $Q$ as $\{V^g_i\}$ depends continuously on
$\mathcal{P}$ by Theorem A3; $Q^N$ is a compact invariant set;
${\dot V}$ is negative semi-definite in $Q^N$; $E = \dot{V}^{-1}(0)
= \{\tilde{C}_{V^g_{i}}\}$, which itself is the largest invariant
subset of $E$ by the control law (\ref{ctrl1_HLOP}). Thus, by
LaSalle's invariance principle, the trajectories of the agents
governed by control law (\ref{ctrl1_HLOP}), starting from any
initial configuration $\mathcal{P}(0) $ (Note that $p_i(0) \in Q,
\forall i \in I_N$), will asymptotically converge to set $E$, the
critical points of $\mathcal{H}$, that is, the centroids of the
generalized centroidal Voronoi cells with respect to the density
function as perceived by the sensors. \hfill $\Box$

Note that a similar result is provided in \cite{bullo1,Bullo09} for
homogeneous sensors using standard Voronoi partition. Further, We
had assumed $p_i(0) \neq p_j(0), \forall i \neq j$. If $f_i \neq
f_j, \forall i \neq j$, then it is possible that $p_i(t) = p_j(t)$
at some time $t >0$. However, if for any $i \neq j$, $f_i(\cdot) =
f_j(\cdot)$, then $p_i(t) \neq p_j(t), \forall t>0$. Thus, the state
space is not the entire $Q^N$, and hence not necessarily compact. An
extension to LaSalle's invariance principle in such a situation is
provided in \cite{Bullo09}, which can be used here for proving the
convergence result, in such a scenario.

It should be noted that the gradient based control law
(\ref{ctrl1_HLOP}) can only guarantee local maximum of the objective
function (\ref{HLOP}). Techniques such as deterministic annealing
can be used to solve the global optimization problem
\cite{salapaka,Rose98,Sharma08}.

\section{Limited range sensors}
In the previous sections, it was assumed that the sensors have
infinite range but with diminishing effectiveness. However, in
reality the sensors will have limited range. In this section a
spatially distributed limited range locational optimization problem
is formulated. In \cite{bullo2} authors address the effect of limit
on the sensor range using standard Voronoi partition. The results
provided here are based on and extension of those results.

Let $R_i$ be the limit on range of the sensors and
$\bar{B}(p_i,R_i)$ be a closed ball centered at $p_i$ with a radius
of $R_i$. The $i$-th sensor has access to information only from
points in the set $V^g_i \cap \bar{B}(p_i,R_i)$. Consider $\hat f_i$
which is continuously differentiable in $[0,R_i]$, with
$\hat{f}_i(R_i) = 0$, and $\hat f_i(r) = 0, \forall r > R_i$ as the
sensor effectiveness function. This function models the
effectiveness of a sensor having a limit of $R_i$ on its range. Note
that the derivative of $\hat f_i$ can have a discontinuity at $r =
R_i^+$, where $R_i^+ = \lim_{h \rightarrow 0} (R_i+h)$. Consider the
following objective function to be maximized,
\begin{equation}
\label{obj_sensr_rng} \hat{\mathcal{H}}(\mathcal{P}) = \sum_{i \in
I_N}\int_{(V^g_{i}\cap\bar{B}(p_i,R_i))}\phi(q)\hat{f_i}(\| p_i - q
\|))dQ
\end{equation}

It can be noted that the objective function is made up of sums of
the contributions from sets $V^g_i \cap \bar{B}(p_i,R_i)$, enabling
the sensors to solve the optimization problem in a spatially
distributed manner. By generalizing the results in \cite{bullo2}, it
can be shown that the gradient of the multi-center objective
function (\ref{obj_sensr_rng}) with respect to $p_i$ is given by
\begin{equation}
\label{grad_sens_rng} \frac{\partial \hat{\mathcal{H}}}{\partial
p_i} =
\int_{(V^g_{i}\cap\bar{B}(p_i,R_i))}\phi(q)\frac{\partial}{\partial
p_i}\hat{f_i}(\| p_i - q \|)dQ
\end{equation}
Thus, the gradient of the objective function (\ref{obj_sensr_rng})
is
\begin{equation}
\label{grad_sensr_rng} \frac{\partial \hat{\mathcal{H}}(\mathcal{P})}{\partial
p_i} =
-\tilde{M}_{(V^g_{i}\cap\bar{B}(p_i,R_i))}(
p_i-\tilde{C}_{(V^g_{i}\cap\bar{B}(p_i,R_i))})
\end{equation}
where, the mass $\tilde{M}_{(V^g_{i}\cap\bar{B}(p_i,R_i))}$ and the
centroid $\tilde{C}_{(V^g_{i}\cap\bar{B}(p_i,R_i))}$ are now
computed within the region $(V^g_{i}\cap\bar{B}(p_i,R_i))$, that is,
the region of Voronoi cell $V^g_i$, which is accessible to the
$i$-th robot. The critical points/configurations are  $\mathcal{P}^*
= \{\mathcal{P} | p_i =
\tilde{C}_{(V^g_{i}\cap\bar{B}(p_i,R_i))}(\mathcal{P}), i \in
I_N\}$.

Consider the following control law to guide the agents toward the respective centroids
\begin{equation}
\label{ctrl_sensr_rng} u_i = -k_{prop}(p_{i} -
\tilde{C}_{(V^g_{i}\cap\bar{B}(p_i,R_i))})
\end{equation}
%

Further, by generalizing the corresponding results in \cite{bullo2},
it can be shown that the trajectories of the sensors governed by the
control law (\ref{ctrl_sensr_rng}), starting from any initial
condition $\mathcal{P}(0)$, will asymptotically converge to the
critical points of $\hat{\mathcal{H}}$.

\section{Simulation Experiments}
In this section we provide results of a set of simulation
experiments to validate the performance of the proposed deployment
strategy. The space is discretized into hexagonal cells, which are
preferred over square or rectangular cells because of the fact that
in hexagonal grids all the neighboring cells are at equal distance.
We have considered  $f_i(r_i) = k_i e^{-\alpha_i {r_i}^2}$ as node
functions. Three density distributions have been used: i) a uniform
density, ii) an exponential density distribution $\phi(x,y) =
0.9e^{-0.001((x-42)^2 + (y-36)^2)}$ having a peak near the center of
space, and iii) an exponential density distribution $\phi(x,y) =
0.9e^{-0.001((x-70)^2 + (y-60)^2)}$ having peak toward the corner of
the space. The space considered is a rectangular area with x-axis
range of 0 to 86 units and y-axis range of 0 to 74 units. The
(locally) optimal deployment is achieved when each of the agents are
located at the respective centroids. Each hexagonal cell measures 1
unit from center to any of its corners. A normalized error measure
of how close the agents are to the centroid at the $n$-th time step
is defined as $error(n) = (\sum_{m=1}^N
\|(p_m(n)-\tilde{C}_{V^g_m}(n))\|^2)/\max_{i\in I_M}(error(i))$,
where $M$ is the total number of time steps to achieve the optimal
deployment. Simulations were carried out by either fixing $k_i$ or
$\alpha_i$, or varying both. The case of both $k_i$ and $\alpha_i$
fixed corresponds to the homogeneous agent case and leads to a
standard Voronoi partition.

\begin{figure*}
\centerline{
\subfigure[]{\psfig{figure=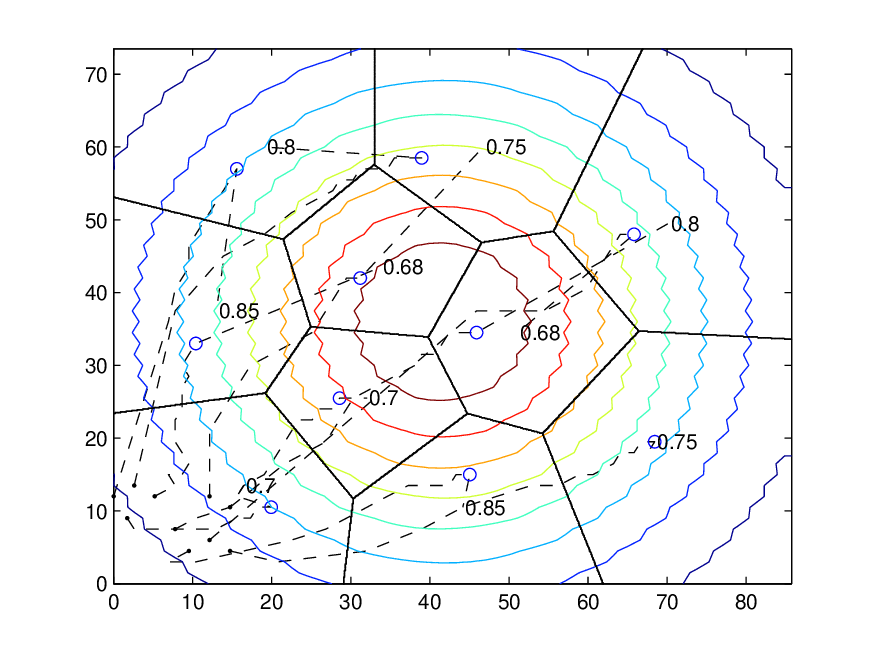,height=5.5cm,width=5.5cm}}
\subfigure[]{\psfig{figure=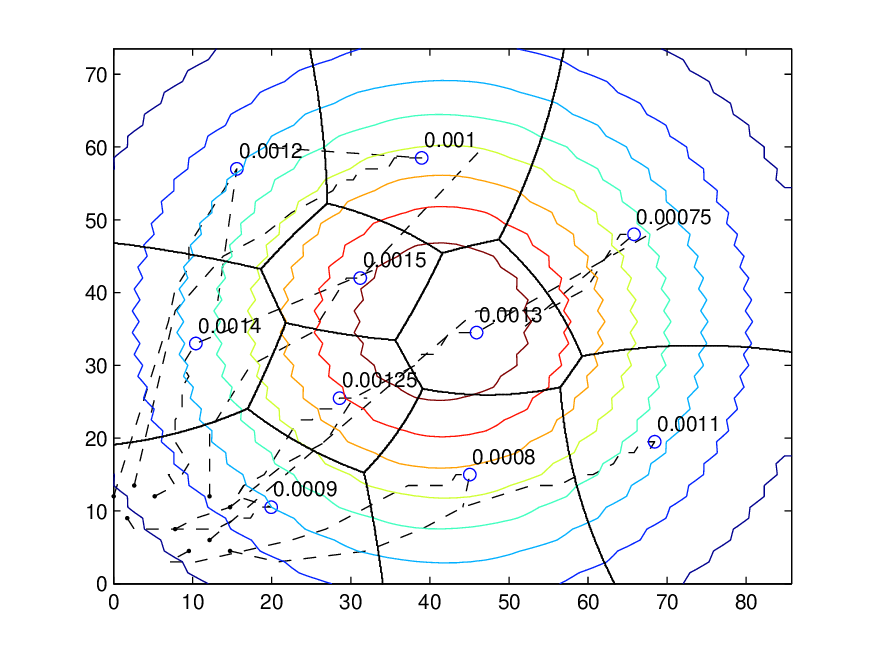,height=5.5cm,width=5.5cm}}
\subfigure[]{\psfig{figure=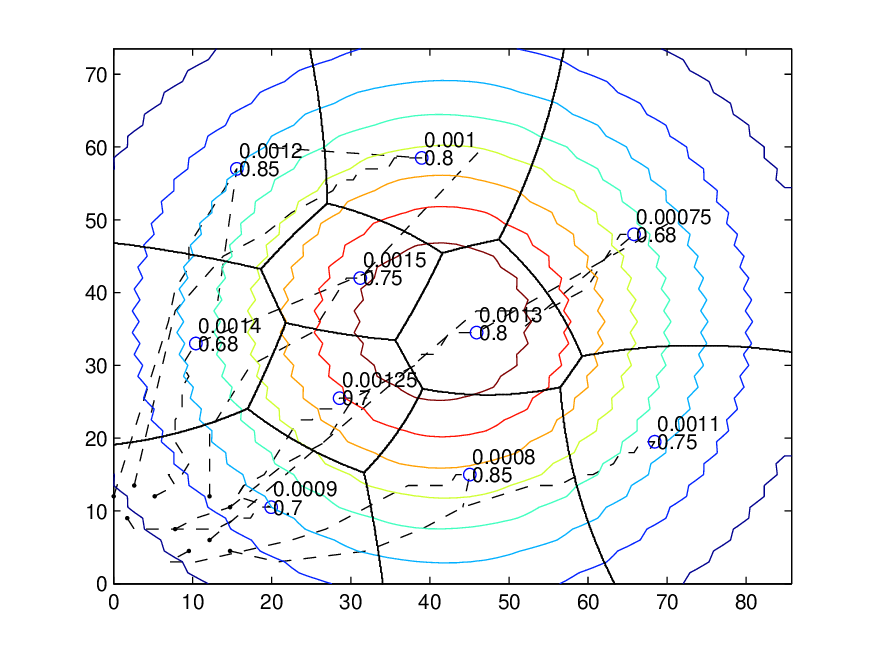,height=5.5cm,width=5.5cm}}
}\caption{Trajectories of agents with exponential density having
peak at the center of space, and a) $\alpha$ fixed and $k$ varying
b) $\alpha$ varying and $k$ fixed c) both $\alpha$ and $k$
varying.}\label{set2_exp_cen}
%
\centerline{
\subfigure[]{\psfig{figure=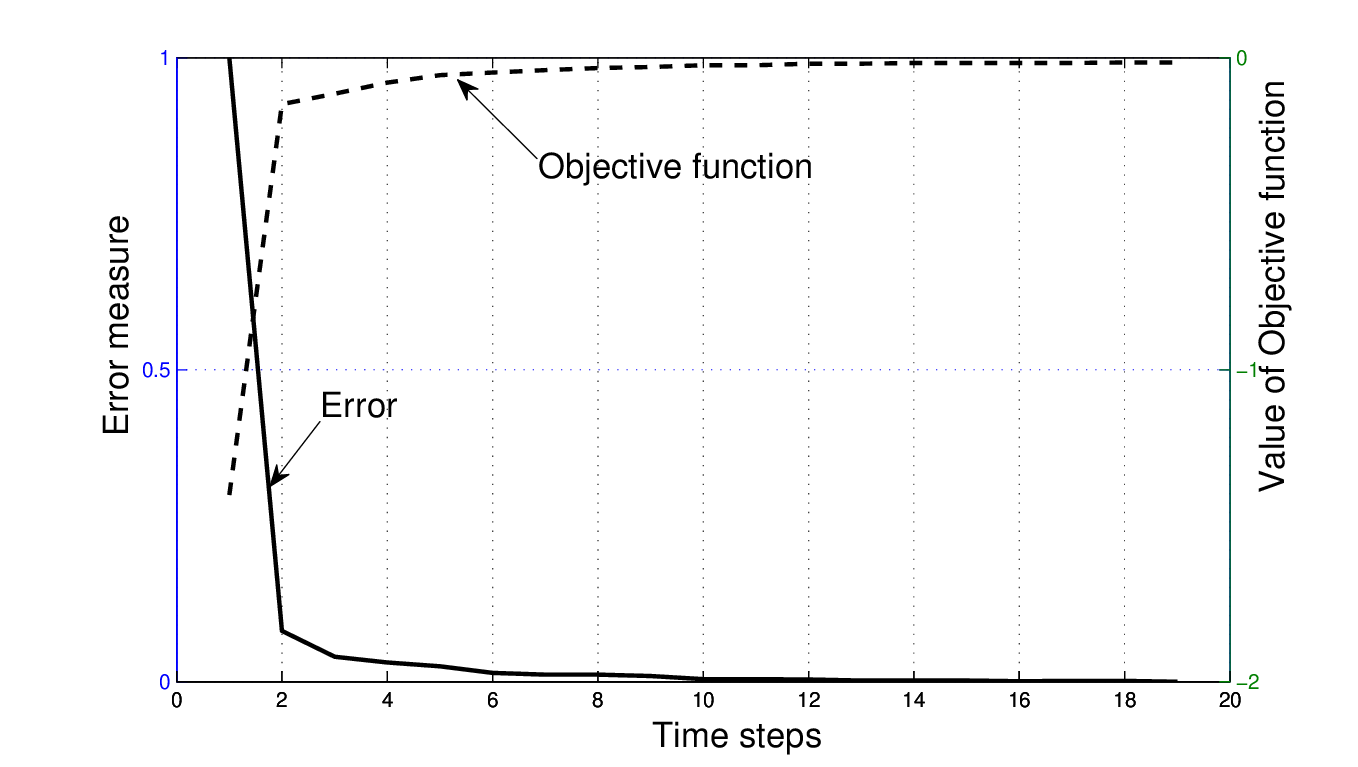,height=5.5cm,width=5.5cm}}
\subfigure[]{\psfig{figure=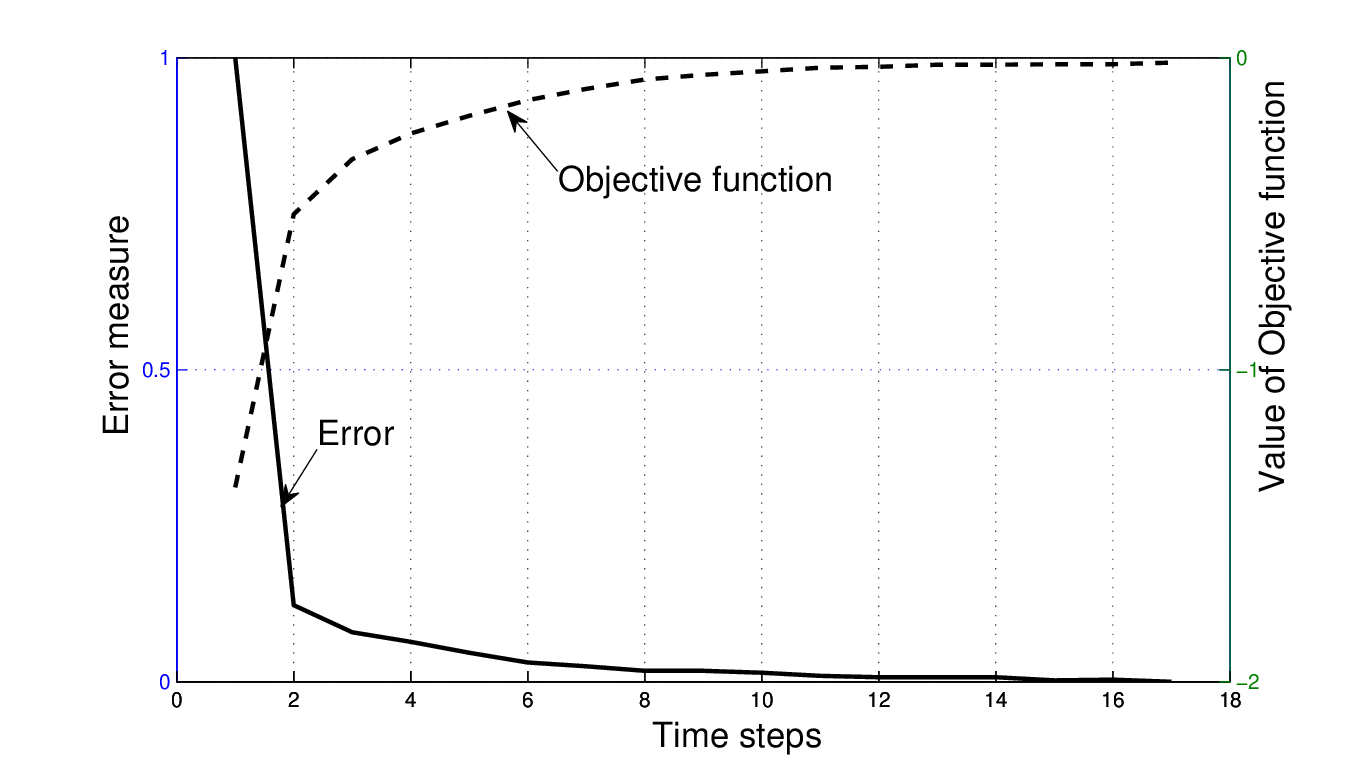,height=5.5cm,width=5.5cm}}
\subfigure[]{\psfig{figure=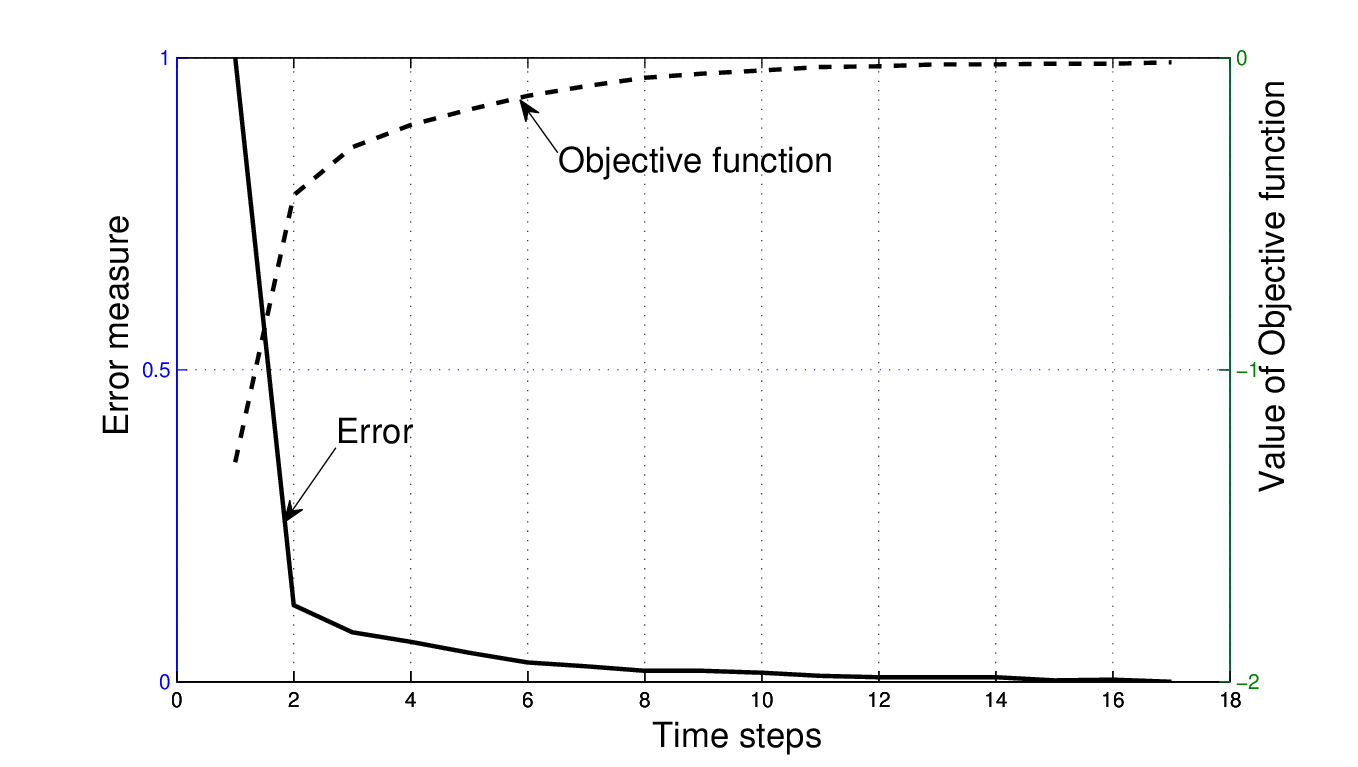,height=5.5cm,width=5.5cm}}
}\caption{History of error and value of objective function with
exponential density having peak at the center of space, and a)
$\alpha$ fixed and $k$ varying b) $\alpha$ varying and $k$ fixed c)
both $\alpha$ and $k$ varying.}\label{set2_exp_cen_err_obj}
\end{figure*}

\begin{figure*}
\centerline{
\subfigure[]{\psfig{figure=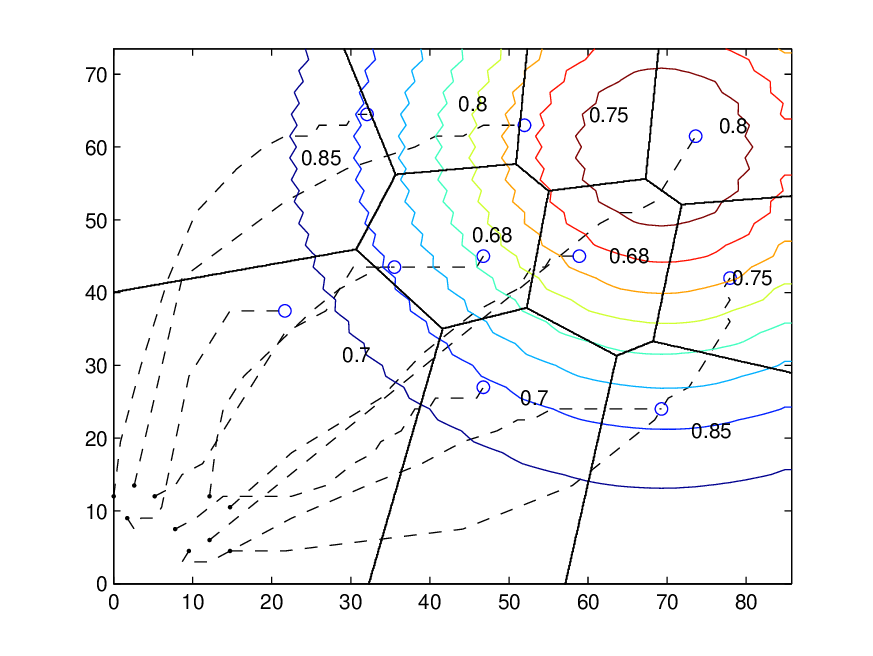,height=5.5cm,width=5.5cm}}
\subfigure[]{\psfig{figure=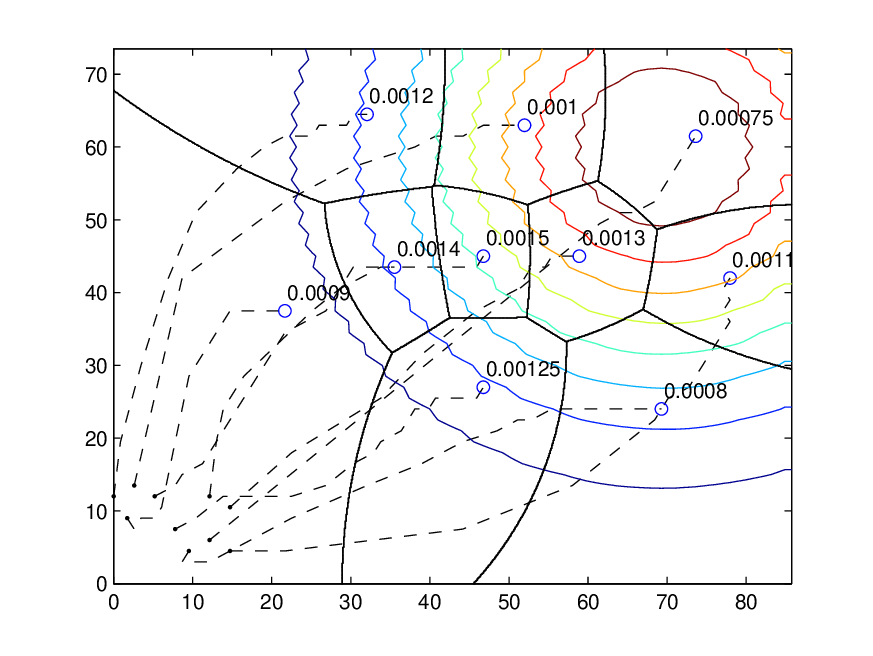,height=5.5cm,width=5.5cm}}
\subfigure[]{\psfig{figure=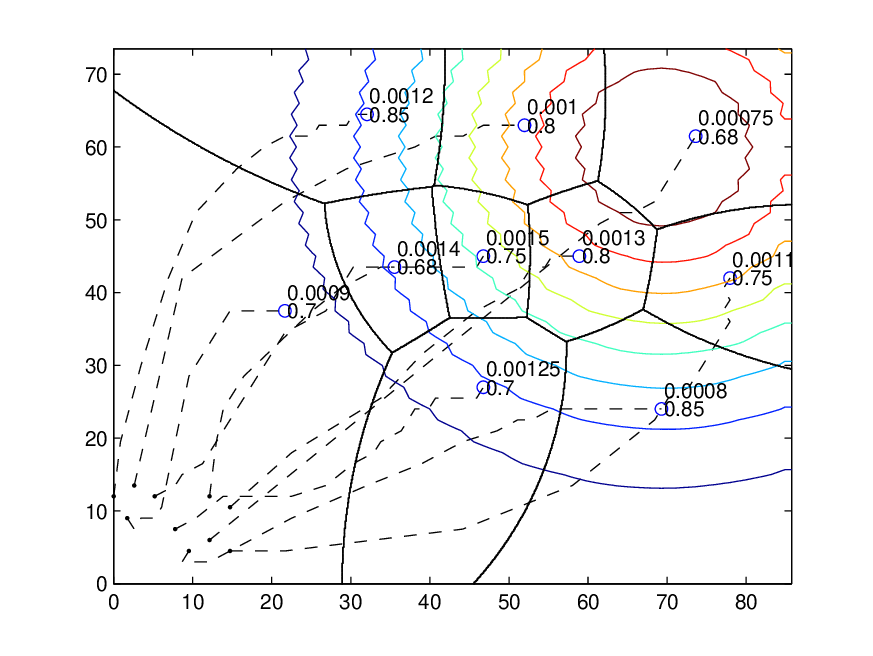,height=5.5cm,width=5.5cm}}
}\caption{Trajectories of agents with exponential density having
peak toward the corner of space, and a) $\alpha$ fixed and $k$
varying b) $\alpha$ varying and $k$ fixed c) both $\alpha$ and $k$
varying.}\label{set2_exp_corner} \centerline{
\subfigure[]{\psfig{figure=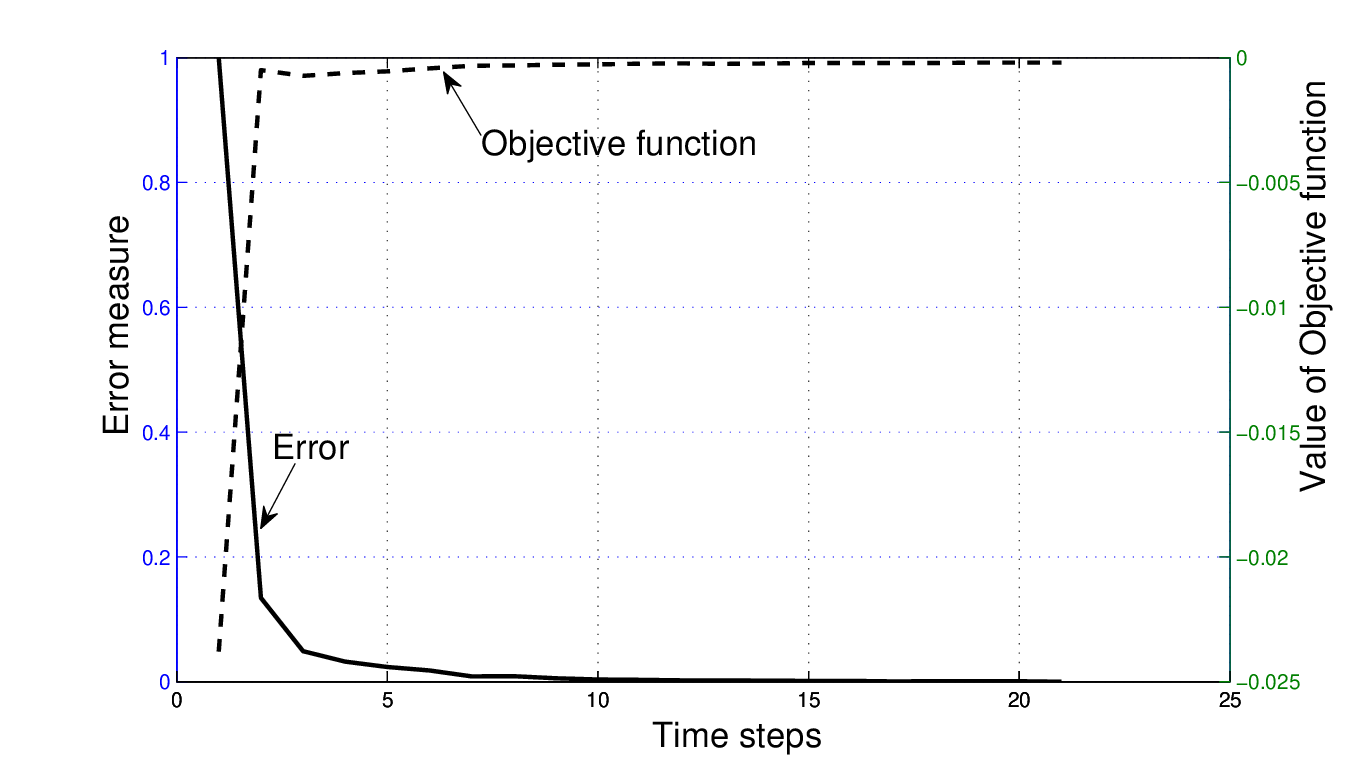,height=5.5cm,width=5.5cm}}
\subfigure[]{\psfig{figure=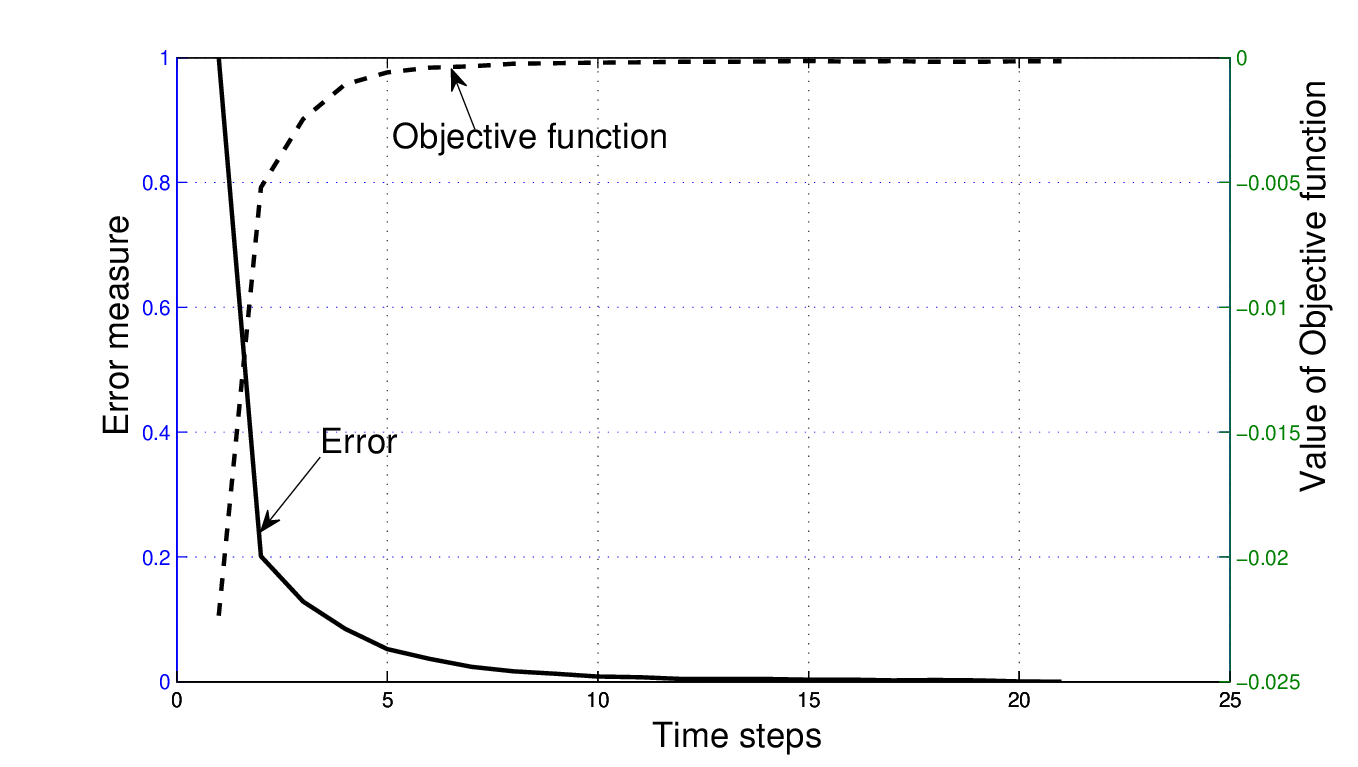,height=5.5cm,width=5.5cm}}
\subfigure[]{\psfig{figure=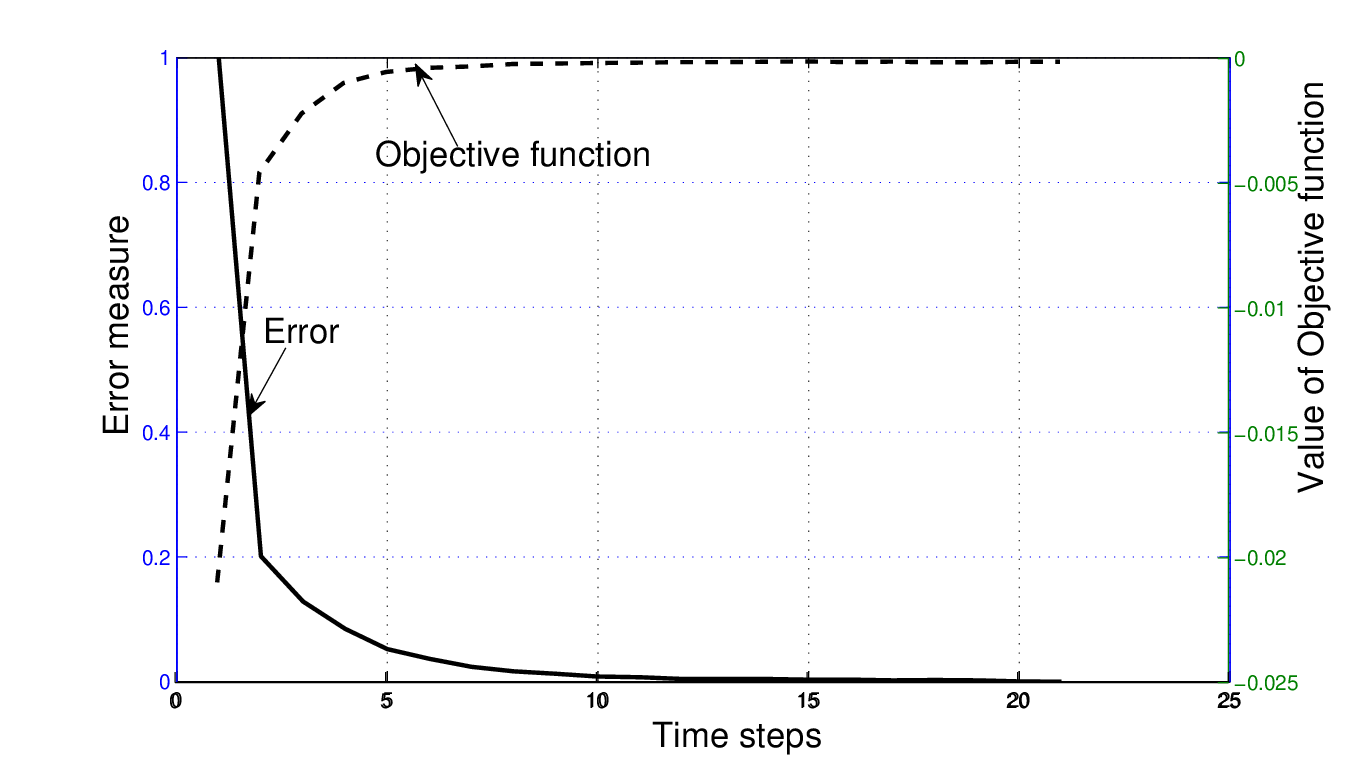,height=5.5cm,width=5.5cm}}
}\caption{History of error and value of objective function with
exponential density having peak toward the corner of space, and a)
$\alpha$ fixed and $k$ varying b) $\alpha$ varying and $k$ fixed c)
both $\alpha$ and $k$ varying.}\label{set2_exp_corner_err_obj}
\end{figure*}

\begin{figure*}
%
\centerline{
\subfigure[]{\psfig{figure=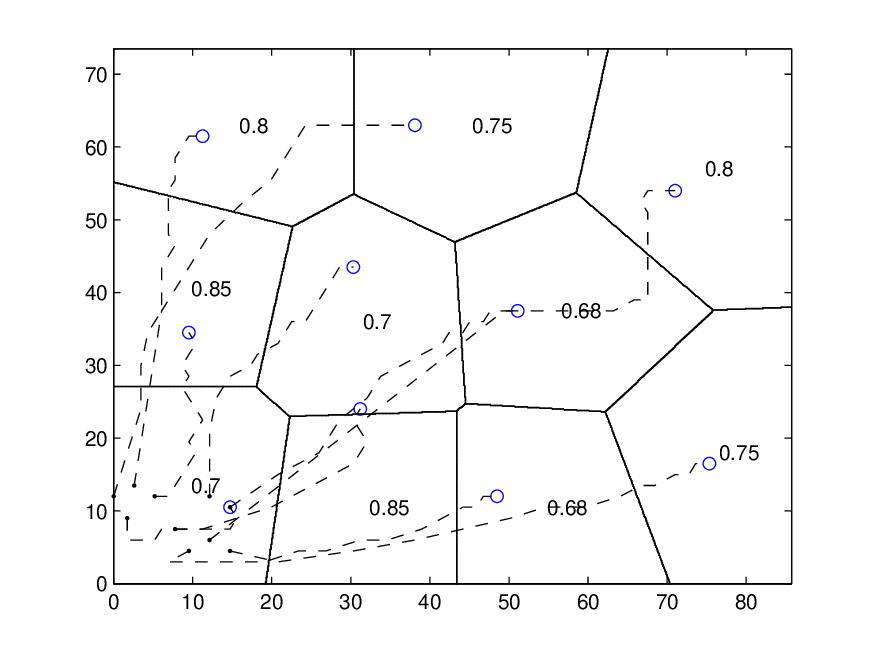,height=5.5cm,width=5.5cm}}
\subfigure[]{\psfig{figure=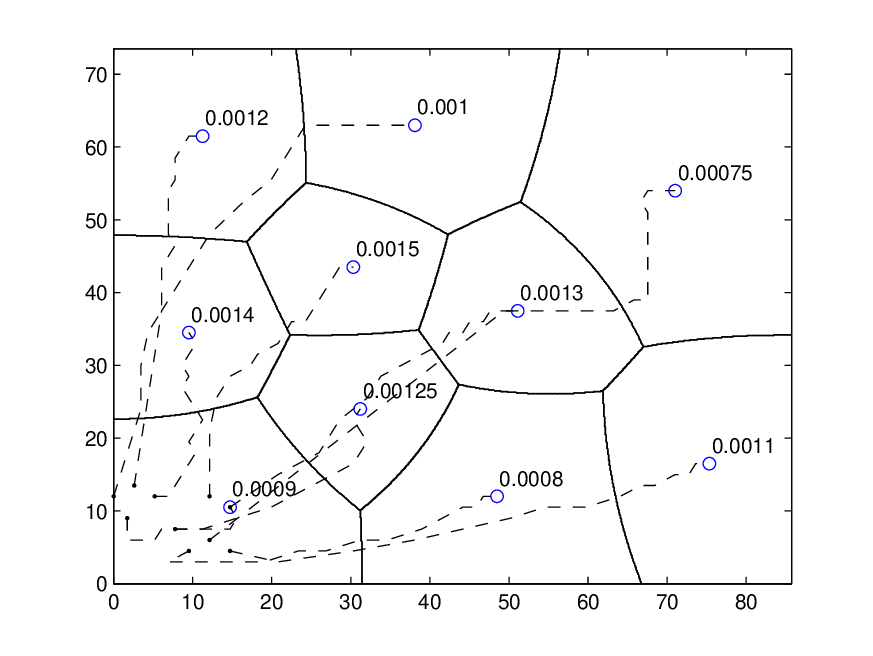,height=5.5cm,width=5.5cm}}
\subfigure[]{\psfig{figure=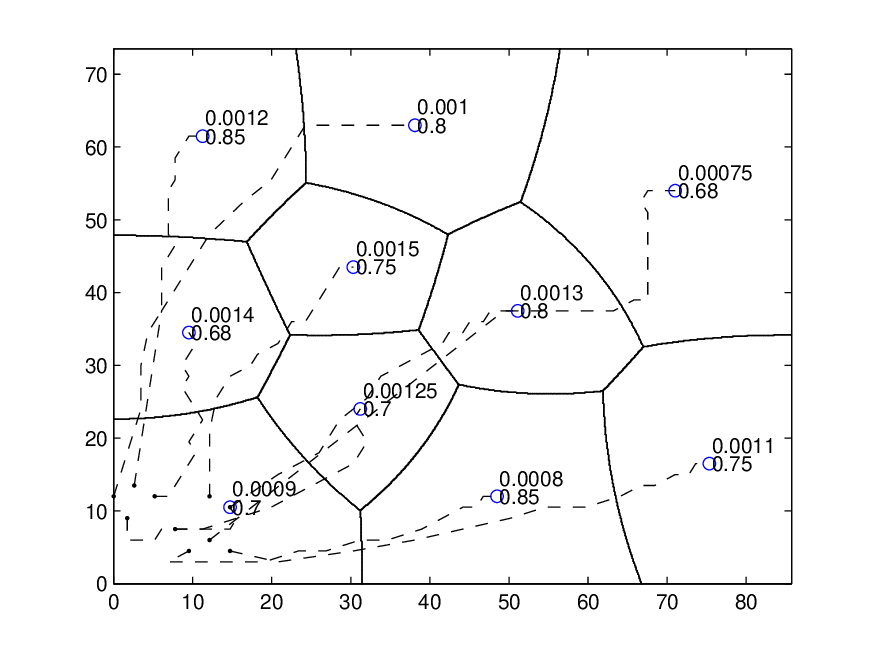,height=5.5cm,width=5.5cm}}
}\caption{Trajectories of agents with uniform density, and a)
$\alpha$ fixed and $k$ varying b) $\alpha$ varying and $k$ fixed c)
both $\alpha$ and $k$ varying.}\label{set2_uniform}
%
\centerline{
\subfigure[]{\psfig{figure=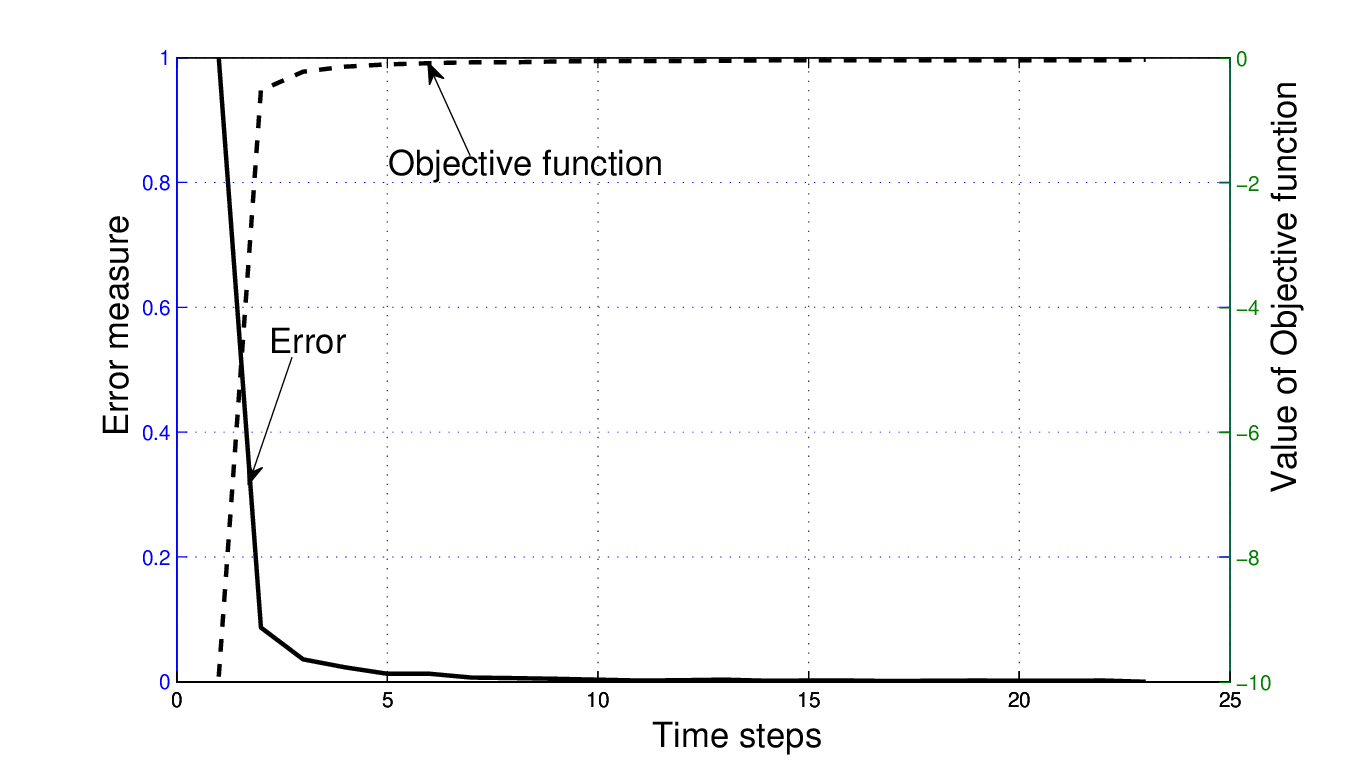,height=5.5cm,width=5.5cm}}
\subfigure[]{\psfig{figure=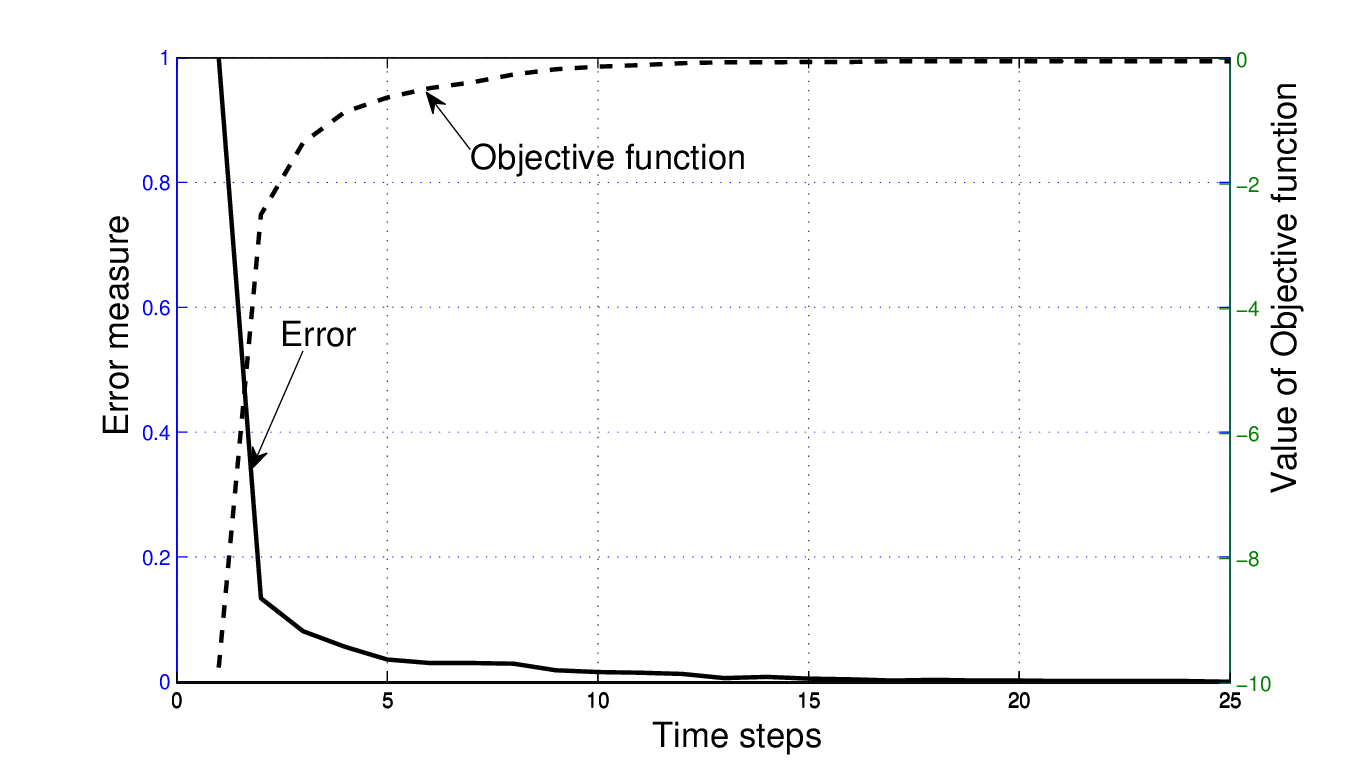,height=5.5cm,width=5.5cm}}
\subfigure[]{\psfig{figure=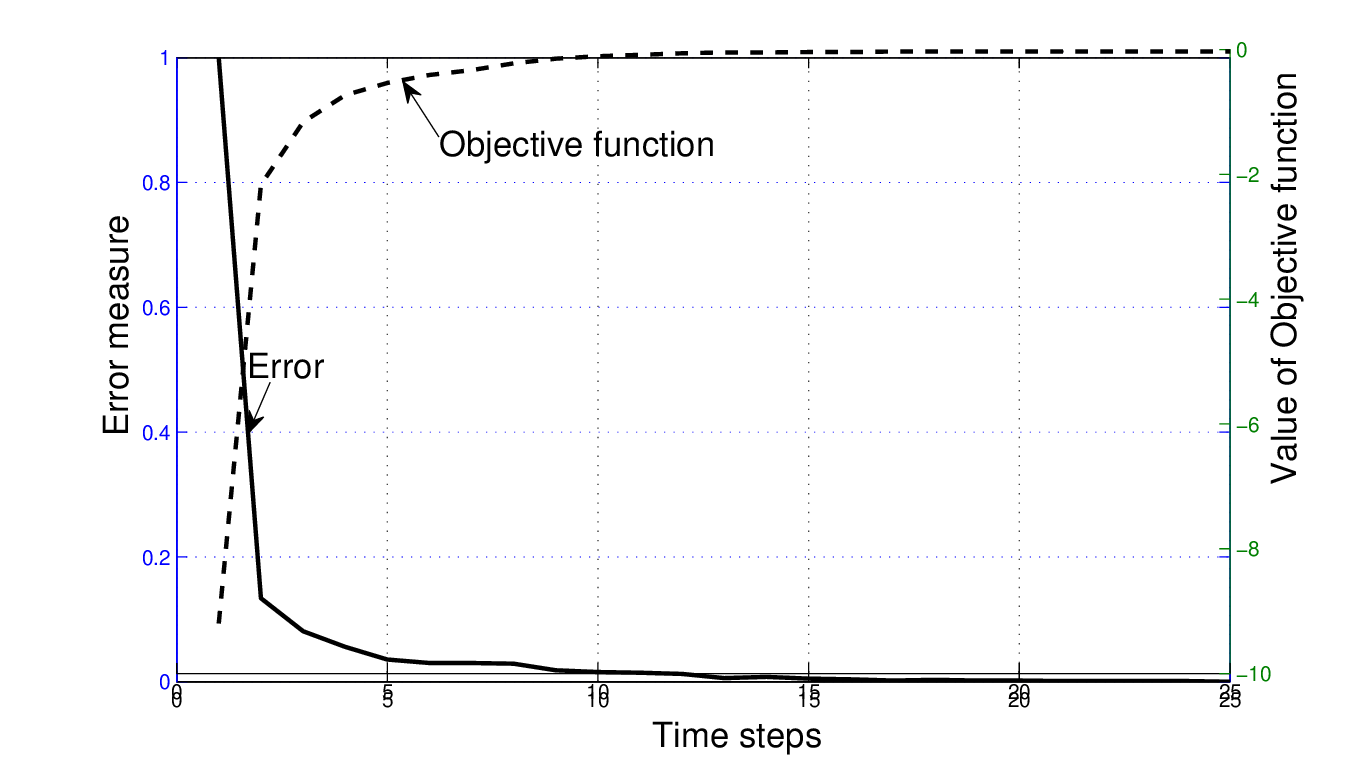,height=5.5cm,width=5.5cm}}
}\caption{History of error and value of objective function with
uniform density, and a) $\alpha$ fixed and $k$ varying b) $\alpha$
varying and $k$ fixed c) both $\alpha$ and $k$
varying.}\label{set2_uniform_err_obj}
\end{figure*}

Figures \ref{set2_exp_cen}, \ref{set2_exp_corner}, and
\ref{set2_uniform} show the trajectories of agents with exponential
density having peak near the center, exponential density having peak
near the corner, and an uniform density. Figures
\ref{set2_exp_cen}(a), \ref{set2_exp_corner}(a), and
\ref{set2_uniform}(a) correspond to cases where all the agents have
the same $\alpha_i$, but different $k_i$. Figures
\ref{set2_exp_cen}(b), \ref{set2_exp_corner}(b), and
\ref{set2_uniform}(b) correspond to case where all the agents have
the same $k_i$ and different $\alpha_i$. Figures
\ref{set2_exp_cen}(c), \ref{set2_exp_corner}(c), and
\ref{set2_uniform}(c) correspond to cases where all the agents have
different $k_i$ and $\alpha_i$. When $k_i$ is fixed, it leads to
multiplicatively weighted Voronoi partition, when $\alpha_i$ is
fixed it leads to a generalized Voronoi partition with intersection
between any two cells being a straight line segment, and when both
$k_i$ and $\alpha_i$ vary, it leads to a generalized Voronoi
partition. A contour plot is also shown for exponential density
cases. It can be observed that the agents move toward the peak of
density. Dots indicate the starting position of agents and `o's
indicate the final position. Figures \ref{set2_exp_cen_err_obj},
\ref{set2_exp_corner_err_obj}, and \ref{set2_uniform_err_obj} show
the history of error and the objective functions for the
corresponding simulations. In all cases the error reduces and the
objective function is (locally) maximized.


\section{Conclusions} A generalization of Voronoi partition has been proposed and the
standard Voronoi decomposition and its variations are shown to be
special cases of this generailzation. The problem of optimal
deployment of agents carrying sensors with heterogeneous
capabilities has been formulated and solved using the generalized
Voronoi partition. The generalized centroidal Voronoi configuration
was shown to be a locally optimal configuration. Illustrative
simulation results were provided to support the theoretical results.

The generalization of Voronoi partition and the heterogeneous
locational optimization techniques can be applied to a wide variety
of problems, such as spraying insecticides, painting by multiple
robots with heterogeneous sprayers, and many others. The
heterogeneous locational optimization problems can find applications
outside robotics field. One such problem is of deciding on optimal
location for public facilities.
The generalization of Voronoi partition presented also provides new
challenges for developing efficient algorithms for computations
related to the generalized Voronoi partition, and characterizing its
properties.

\section*{Acknowledgements}
This work was partly supported by Indian National Academy of
Engineering (INAE) fellowship under teacher mentoring program.

\section*{Appendix: On continuity of generalized Voronoi partition}
A proof on the continuity of the generalized Voronoi partition is
provided here. The treatment provided here is kept informal, but
discusses major steps for a more elaborate mathematical proof, which
is beyond the scope of this paper.

Consider a partition $\mathcal{W}(\mathcal{P}) =
\{W_i(\mathcal{P})\}\text{,} i \in I_N$ of $Q \subset \mathbb{R}^2$
with $\mathcal{P} = \{p_1,p_2,\ldots,p_N \}$, $p_i \in Q$. Let
$\partial{W}_i(\mathcal{P})$ be the bounding (closed) curve of cell
$W_i$, parameterized by a parameter $u_i(\mathcal{P}) \in [0,1]$.
The partition $\mathcal{W}(\mathcal{P})$ is said to be continuous in
$\mathcal{P}$ if and only if each of the $u_i(\mathcal{P})$ is
continuous in $\mathcal{P}$. That is, the cell boundaries move
smoothly with configuration of sites $\mathcal{P}$. This is an
important property needed for proving the stability and convergence
of the agent motion under the influence of the proposed control law
using LaSalle's invariance principle. This understanding of
continuity of a partition can be extended to a general
$d$-dimensional Euclidean space. Now, in $\mathbb{R}^d$, $\partial
W_i$ are hypersurfaces in $d$-dimensional Euclidean space and $u_i$
are vector valued functions of dimension $d-1$. Continuity of
$u^j_i$, the $j$-th component of $u_i$ with $\mathcal{P}$ for each
$i \in I_N$ and $j\in I_{d-1}$ ensures the continuity of
$\mathcal{W}(\mathcal{P})$ when $Q \subset \mathbb{R}^d$.

\begin{figure}
\centerline{
\subfigure[]{\psfig{figure=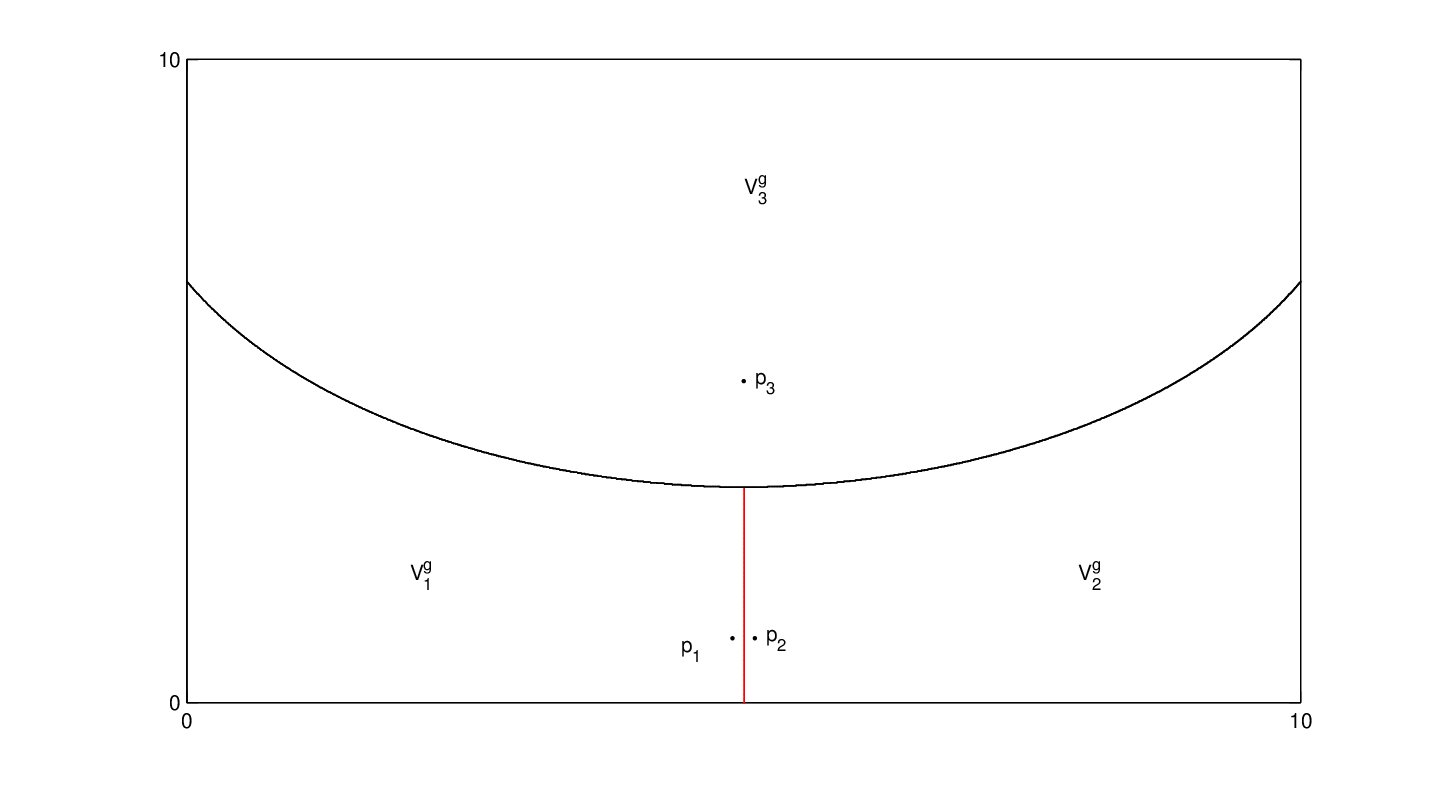,height=5.5cm,width=5.5cm}}
\subfigure[]{\psfig{figure=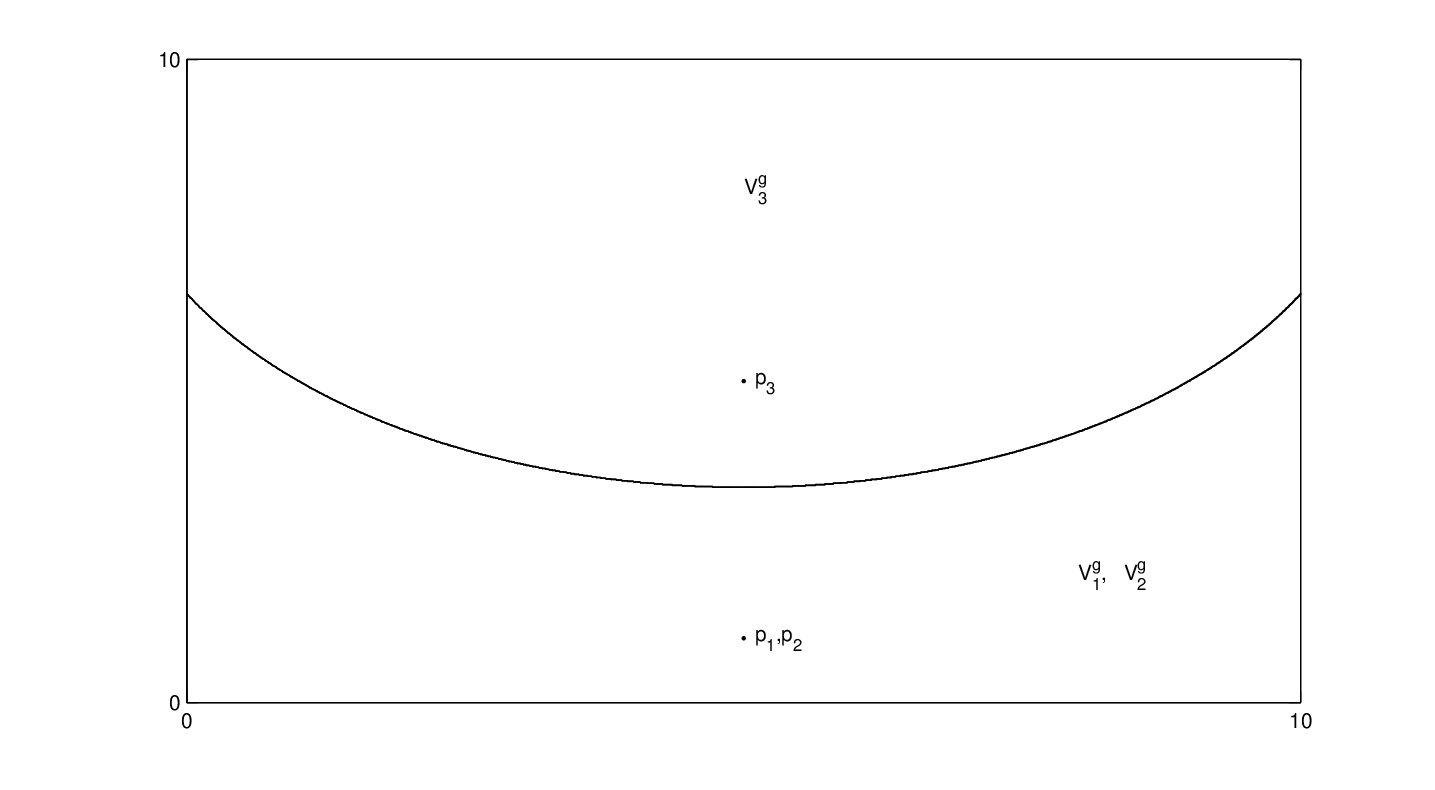,height=5.5cm,width=5.5cm}}
\vspace{-0.1in} }\caption{Illustration of discontinuity of Voronoi
partition at $p_1 = p_2$ when $f_1 = f_2$. The straight line segment
separating $V^g_1$ and $V^g_2$ suddenly disappears when two points
$p_1$ and $p_2$ merge.}\label{vor_cont_illustrate}
\end{figure}

Now consider the continuity of standard Voronoi partition. It is
well known that as long as $p_i \neq p_j$, whenever $i \neq j$,
$\mathcal{V}(\mathcal{P})$, the standard Voronoi partition, is
continuous in $\mathcal{P}$. This continuity is due to continuity of
the Euclidean distance. It is interesting to see what happens when
$p_i = p_j$ for some $i \neq j$. Consider $i \neq j$ and let $p_i(0)
\neq p_j(0)$. Now, $V_i$ and $V_j$ are distinct sets such that $V_i
\cap V_j$ is either null or is the common boundary between them,
which is a straight line segment of the perpendicular bisector of
the line joining $p_i$ and $p_j$.  Let $V_i \cap V_j \neq
\emptyset$, that is, the agents $i$ and $j$ are Voronoi neighbors,
and agent $i$ is stationary while  agent $j$ is moving toward $i$.
This causes $V_i$ to shrink and $V_j$ to expand as the common
boundary between them moves closer to $p_i$. But, when $p_j(t) =
p_i(t)$, $V_i(t) = V_j(t) = V_i(0) \cup V_j(0)$, which is a sudden
jump, and boundaries of both $V_i$ and $V_j$ jump, leading to
discontinuity. Thus, $\mathcal{V}(\mathcal{P})$ is discontinuous
whenever there is a transition between $p_i \neq p_j$ and $p_i =
p_j$, for some $i \neq j$. Similar discontinuity, as in the standard
Voronoi partition with $\mathcal{P}$, can occur in the case of
generalized Voronoi partition, whenever $p_i$ and $p_j$ are Voronoi
neighbors and $f_i = f_j$. This is illustrated in Figure
\ref{vor_cont_illustrate}. The line segment separating $V^g_1$ and
$V^g_2$ disappears when $p_1 = p_2$ and the corresponding Voronoi
cells jump in their size. In fact $V^g_1$ and $V^g_2$ jump to $V^g_1
\cap V^g_2$, and the common boundary between $V^g_1$ and $V^g_2$
jumps to the boundary of $V^g_1 \cap V^g_2$.

In case of the standard Voronoi partition, if $p_i(0) \neq p_j(0)$,
whenever $i \neq j$, the control law ensures that $p_i(t) \neq
p_j(t)$, whenever $i \neq j$, $\forall t > 0$. If two agents are
together to start with, that is, $p_i(0) = p_j(0)$, whenever $i \neq
j$, then the control law ensures that they are always together. This
is due to the fact that $\tilde{C}_{V_i} \in V_i$, $\forall i \in
I_N$. This is not always true in case of the generalized Voronoi
partition. Thus, the issue of continuity of generalized Voronoi
partition is a more involved one.


\nnd {\bf Condition A.} There is no transition from $\mathcal{P}(t)$
to $\mathcal{P}(t')$, such that $p_i(t) \neq p_j(t)$  and $p_i(t') =
p_j(t')$, for all pairs $(i,j)$, $i \neq j$, $i,j \in I_N$, and $f_i
= f_j$.

As discussed earlier, violation of the condition A can cause discontinuity in the generalized Voronoi partition.

\nnd {\it Lemma A1:~}The generalized Voronoi partition depends continuously on $\mathcal{P}$ if condition A is satisfied.

\nnd {\it Outline of the proof:~} Let us consider any two agents $i$
and $j$ which are Voronoi neighbors, take up three distinct and
exhaustive cases, and prove the continuity for each of them.

\noindent {\bf Case i)} No two node functions are identical, that
is, $(f_i-f_j) \neq 0$, whenever $i \neq j$, and $V^g_i \neq
\emptyset$, $\forall i \in I_N$. $\forall \mathcal{P}$.

All the points $q \in Q$ on the boundary common to $V^g_i$ and
$V^g_j$ is given by  $\{q \in Q | f_i(\|p_i-q\|)=f_j(\|p_j-q\|)\}$,
that is, the intersection of the corresponding node functions. Let
the $j$-th agent move by a small distance $dp$. This makes a point
$q \in Q$, on the common boundary between $V_i$ and $V_j$ move by a
distance, say $dx$. Now, as the node functions are strictly
decreasing and are continuous, it is easy to see that $dx
\rightarrow 0$ as $dp \rightarrow 0$. This is true for any two $i$
and $j$, and any $q$ on the common boundary between $V^g_i$ and
$V^g_j$. Thus, the generalized Voronoi partition depends
continuously on $\mathcal{P}$. Note that it is easy to verify that
$p_i = p_j$ for $i \neq j$ does not lead to discontinuity in this
case.

\nnd {\bf Case ii)} No two node functions are identical, that is,
$(f_i-f_j) \neq 0$, whenever $i \neq j$,  and $\mathcal{P}$ is such
that $V^g_i = \emptyset$ for some $i \in I_N$.

Let both $V^g_i$ and $V^g_j$ be non null sets, and $V^g_i \subset
V^g_j$ for some $\mathcal{P}$. Now with evolution of $\mathcal{P}$
in time, let $V^g_i = \emptyset$ for some $\mathcal{P'} =
\mathcal{P} + \delta \mathcal{P}$. The quantity $\int_{V^g_i} dQ$,
and the common boundary vanish gradually to zero as the agent
configuration changes from $\mathcal{P}$ to $\mathcal{P'}$. Further,
there is no jump in the boundary or size of the (generalized)
Voronoi cell, or in the common boundary  as in the case of standard
Voronoi cell when there is a transition from $p_i \neq p_j$ to $p_i
= p_j$.

\nnd {\bf Case iii)}. For some $i \neq j$, $f_i = f_j$. In this
case, the common boundary between $V^g_i$ and $V^g_j$ is a segment
of the perpendicular bisector of line joining $p_i$ and $p_j$. As in
the case of standard Voronoi diagram, if agents $i$ and $j$ are
neighbors and there is a transition from $p_i \neq p_j$ to $p_i =
p_j$, the $\mathcal{V}^g(\mathcal{P})$ will have discontinuity.
However, if condition A is satisfied, such a discontinuity will not
occur.

Thus, as long as the condition A is satisfied, the generalized
Voronoi partition depends continuously on $\mathcal{P}$. \hfill
$\Box$


\nnd {\it Lemma A2:~} The control law (\ref{ctrl1_HLOP}) ensures
that condition A is satisfied for all $\mathcal{P}$ if $p_i(0) \neq
p_j(0)$ for all pairs $(i,j), i,j \in I_N, i \neq j$, for which $f_i
= f_j$.

\nnd {\it Proof.} Condition A can be violated only for pair $(i,j)$
for which $f_i = f_j$. As $p_i(0) \neq p_j(0)$, condition A will be
violated if at some time $t$, $p_i(t) = p_j(t)$. Now
$\tilde{C}_{V^g_i}$ will lie within the half plane $\{q \in Q |
\|q-p_i \| < \|q-p_j\|\}$ and hence the control law cannot make the
agent $i$ cross the common boundary between $i$ and $j$. Similarly,
agent $j$ too cannot cross the boundary. Hence, if $p_i(0) \neq
p_j(0)$, then $p_i(t) \neq p_j(t)$ for any $t \in \mathbb{R}$. Thus,
the condition A cannot be violated. \hfill$\Box$


\nnd {\it Theorem A3:~} If $\mathcal{P}(0)$ is such that $p_i(0)
\neq p_j(0)$ for all pairs $(i,j), i,j \in I_N, i \neq j$, for which
$f_i = f_j$, and the agents move according to the control law
(\ref{ctrl1_HLOP}), then the generalized Voronoi partition depends
continuously on $\mathcal{P}$.

\noindent{\it Proof.~} The proof follows from Lemmas A1 and A2.
\hfill$\Box$

\end{document}